\documentclass[fleqn]{amsart}
\usepackage{amsmath,amsfonts,amssymb}
\usepackage[a4paper]{geometry}
\def\R{\mathbb R}
\def\C{\mathbb C}
\def\D{\mathrm D}
\def\d{\mathrm d}
\def\i{\mathrm i}
\def\e{\mathrm e}
\DeclareMathOperator{\diag}{diag}
\DeclareMathOperator{\spec}{spec}
\DeclareMathOperator{\trace}{tr}
\newtheorem{thm}{Theorem}[section]
\numberwithin{equation}{section}
\begin{document}
\author{Jens Wirth}
\title[Periodic dissipation]
{On the influence of time-periodic dissipation on energy and dispersive estimates}
\address{Jens Wirth, Department of Mathematics, University College London, Gower Street, London WC1E 6BT, UK}
\email{jwirth@math.ucl.ac.uk}
\begin{abstract}In this short note we discuss the influence of a time-periodic dissipation term on a-priori
estimates for solutions to dissipative wave equations. The approach is based on a diagonalisation argument for high frequencies and results from spectral theory of periodic differential equations / Floquet theory for bounded frequencies. \end{abstract}
\maketitle

In recent years much attention was paid to hyperbolic equations with time-dependent coefficients and
the question of the influence of the precise time-dependence on asymptotic properties of solutions. In particular, we refer to \cite{ReissigYagdjian}, \cite{ReissigHirosawa}, \cite{ReissigSmith}, \cite{Wirth06}, \cite{Wirth07} or the survey article \cite{Reissig} for an overview of results. All these papers have in common that they use assumptions on derivatives of the coefficients to avoid the (bad) influence of oscillations. That oscillations may have deteriorating influences was shown in \cite{Yagdjian01} for the example of a wave equation with time-periodic speed of propagation. In this case (some) solutions have exponentially growing energy. The counter-example of \cite{ReissigSmith} shows that even for the Cauchy problem
$$ u_{tt}-a^2(t)\Delta u = 0 $$
with $a(t)=2+\sin(\log(e+t)^\alpha)$ and $\alpha>2$ (some) solutions have supra-polynomially increasing energy, while for $\alpha=1$ polynomial growth may occur and for $\alpha<1$ the energy can be bounded by $t^{\epsilon}$ for any $\epsilon>0$. Similar results can be obtained for oscillations in mass terms, especially for the case of periodic mass terms we can always find solutions with exponentially increasing energy.

Our aim is to show that the influence of oscillations in the dissipation term is different. To be precise, we
want to show that a wave equation with a periodic in time dissipation term,
$$ u_{tt}-\Delta u+2b(t) u_t=0,\qquad b(t+T)=b(t)\ge0, $$
satisfies the same Matsumura-type estimate as obtained for constant dissipation in \cite{Matsumura76}. 
We conjecture that the results of \cite{Wirth07} (where only 
very slow oscillations were treated and decay results of the same structure were obtained) can be extented to general dissipation terms with $tb(t)\to\infty$ without further assumptions on derivatives.
However, it is an open problem how to achieve such a result.

This note is organised as follows: In Section~\ref{sec1} we give the basic assumptions on the
Cauchy problem under consideration and discuss properties of its fundamental solution and the associated monodromy operator. Main result is Theorem~\ref{thm1} on page~\pageref{thm1}.
In Section~\ref{sec2} we discuss applications of Theorem~\ref{thm1} to energy and 
more generally $L^p$--$L^q$ decay estimates of solutions. Results are given in Theorems~\ref{thm2}
and \ref{thm3}. Finally Theorem~\ref{thm4} implies a diffusion phenomenon for periodically damped 
wave equations.

\section{Representation of solutions}\label{sec1}

We consider the Cauchy problem
\begin{equation}\label{eq:CP}
  u_{tt}-\Delta u+2b(t) u_t =0,\qquad u(0,\cdot)=u_1,\quad u_t(0,\cdot)=u_2
\end{equation}
for a wave equation with time-dependent dissipation, where we assume that the coefficient
$b(t)$ is continuous, of bounded variation and satisfies $b(t)>0$ a.e. together with
\begin{equation}
  b(t+T)=b(t)\quad\text{for some $T>0$ and all $t\in\R$}.
\end{equation}
We denote the mean value of $b(t)$ as
\begin{equation}
  \beta = \frac1T \int_0^T b(t)\d t.
\end{equation}
Using a partial Fourier transform with respect to the spatial variables we reduce the Cauchy problem
to the ordinary differential equation 
\begin{equation}\label{eq:ODE}
\hat u_{tt}+|\xi|^2\hat u+2b(t)\hat u_t=0,
\end{equation}
which we reformulate as first order system $\D_t V = A(t,\xi) V$ for $V=(|\xi|\hat u,\D_t\hat u)^T$ with coefficient matrix
\begin{equation}
  A(t,\xi)=\begin{pmatrix} &|\xi|\\|\xi|&2\i b(t)\end{pmatrix}
\end{equation}
 where $\D_t=-\i\partial_t$ denotes the Fourier derivative. Aim of our investigation is to describe the
 corresponding fundamental solution $\mathcal E(t,s,\xi)$, i.e. the matrix-valued solution to
 \begin{equation}
   \D_t \mathcal E(t,s,\xi) = A(t,\xi) \mathcal E(t,s,\xi),\qquad \mathcal E(s,s,\xi)=I\in\C^{2\times 2}
 \end{equation}
 or, exploiting the periodic structure of the problem, the corresponding family of monodromy
 matrices $\mathcal M(t,\xi)=\mathcal E(t+T,t,\xi)$ for $t\in[0,T]$. Note, that these matrices satisfy
 the periodic problem
 \begin{equation}
   \D_t\mathcal M(t,\xi) = [A(t,\xi),\mathcal M(t,\xi)], \qquad \mathcal M(T,\xi)=\mathcal M(0,\xi) .
 \end{equation}
 
 \subsection{A diagonalisation based approach for high frequencies} We consider large frequencies $|\xi|\ge N$,
 the constant $N$ will be determined later on. In this case the influence of $|\xi|$ in $A(t,\xi)$ seems
 to be stronger than the influence of the (comparatively small) coefficient $b(t)$.  To use this, we apply two transformations to the system. For the first step we use the unitary matrices
 \begin{equation}
  M=\frac1{\sqrt2}\begin{pmatrix}1&-1\\1&1\end{pmatrix},\qquad
  M^{-1} = \frac1{\sqrt2} \begin{pmatrix}1&1\\-1&1\end{pmatrix} = M^*
\end{equation}
and set $V^{(0)} = M^{-1} V$, such that
\begin{equation}
  \D_t V^{(0)} = \left( \begin{pmatrix}|\xi|&\\&-|\xi|\end{pmatrix} + \i b(t) \begin{pmatrix}1&1\\1&1\end{pmatrix}\right) V^{(0)}.
\end{equation}
We denote the first (diagonal) matrix as $\mathcal D(\xi)$ and the remainder term as $R_0(t,\xi)$. For
convenience we set $\mathcal D_1 = \mathcal D+\diag R_0$ and $R_1=R_0-\diag R_0$.
In the next step we follow an idea from \cite [Chapter 1.6]{Eastham} and construct a matrix $N_1(t,\xi)$ subject to
\begin{equation}\label{eq:N1-sys}
  \D_t N_1 = [\mathcal D_1, N_1 ] + R_1 
\end{equation}
and $N_1(0,\xi)=I$ (we could use any starting point $N_1(s,\xi)=I$ here without changing much of the calculation). The choice of $N_1(t,\xi)$ implies that
\begin{equation}
   (\D_t-\mathcal D_1-R_1)N_1 - N_1(\D_t - \mathcal D_1) 
   = \D_t N_1 - [\mathcal D_1,N_1] -R_1 N_1   = R_1 (I-N_1),
\end{equation}
such that with $R_2=-N_1^{-1}R_1(I-N_1)$ the operator equation $(\D_t-\mathcal D_1-R_1)N_1
=N_1(\D_t-\mathcal D_1-R_2)$ holds true. It remains to understand in which sense the remainder $R_2$ is better than the remainder $R_1$ and that it is indeed possible to choose the zone constant $N$ large enough to guarantee the invertibility of $N_1$.

For this we solve \eqref{eq:N1-sys}. Since $\mathcal D_1$ is diagonal, we see that $\D_t\diag N_1 = 0$
and therefore $\diag N_1=I$. The two off-diagonal entries of $N_1=(\begin{smallmatrix}1&n^-\\n^+&1\end{smallmatrix})$ satisfy
\begin{equation}
  \D_t n^\pm(t,\xi) = \pm 2|\xi| n^\pm(t,\xi)+\i b(t),\qquad n^\pm(0,\xi)=0, 
\end{equation}
such that
\begin{equation}
  n^\pm(t,\xi) = \int_0^t\e^{\pm 2\i|\xi|(t-s)} b(s)\d s = \int_0^t \e^{\pm2\i s|\xi|} b(t-s)\d s,
\end{equation}
especially we see that $n^+=\overline{n^-}$ and as Fourier transforms of $b(t-s)1_{[0,t]}(s)$ the
Riemann-Lebesgue lemma implies $n^\pm(t,\xi)\to0$ as $|\xi|\to\infty$ for any {\em fixed} $t$. We show
that this is true uniformly in $t\in[0,2T]$, provided that $b$ is of bounded variation. Indeed, integration 
by parts yields
\begin{align}
   n^\pm(t,\xi) = \frac1{\pm2\i|\xi|} \e^{\pm2\i|\xi|s}b(t-s)\big|_{s=0+}^{t-} + \frac1{\pm2\i|\xi|}\int_0^t 
   \e^{\pm2\i|\xi|s}b'(t-s)\d s \le C (1+t) |\xi|^{-1}
\end{align}
and the assertion follows.
 
Now we are in a position to show that $\|\mathcal M(t,\xi)\|<1$ for $|\xi|>N$ for a sufficiently large constant $N$. Note first, that the transformation matrices satisfy $N_1(t,\xi)\to I$ and therefore $N_1^{-1}(t,\xi)\to I$ uniform in $t\in[0,2T]$ as $|\xi|\to\infty$. Furthermore, the remainder term satisfies
$R_2(t,\xi)\to0$ as $|\xi|\to\infty$ uniformly in $t\in[0,2T]$.  Thus for sufficiently large zone constant $N$ we can achieve that
\begin{equation}
\sup_{t\in[0,T]} \|N_1(t+T)\|\; \exp\big(\int_t^{t+T} \| R_2(s,\xi)\| \d s\big)\; \|N_1^{-1}(t)\| < \exp\beta T
\end{equation}
holds true. We fix this choice of the constant $N$ and construct the fundamental solution $\mathcal E(t,s,\xi)$ to $\D_t-A(t,\xi)$. We start from the transformed version of this equation. The fundamental solution to the diagonal part $\D_t-\mathcal D_1$ is given by 
\begin{equation}
  \frac{\lambda(s)}{\lambda(t)}  \tilde{ \mathcal E}_0(t,s,\xi),
\end{equation}
where $\tilde{\mathcal E}_0(t,s,\xi)=\diag\big( \exp(\i(t-s)|\xi|),\exp(-\i(t-s)|\xi|)\big)$ is related to the free propagator $\mathcal E_0(t,s,\xi) = M\tilde{\mathcal E}_0(t,s,\xi)M^{-1}$ corresponding to $b\equiv0$ and
$\lambda(t) = \exp(\int_0^t b(\tau)\d\tau)$ describes the influence of dissipation. Note that $\lambda(t)\approx \exp\beta t$ as $t\to\infty$. For the fundamental solution to $\D_t-\mathcal D_1-R_2$ we make the {\sl ansatz} $\tilde{\mathcal E}_0(t,s,\xi)\mathcal Q(t,s,\xi)$ such that
\begin{equation}
  \D_t \mathcal Q(t,s,\xi) = \tilde{\mathcal E}_0(s,t,\xi)R_2(t,\xi)\tilde {\mathcal E}_0(t,s,\xi) \mathcal Q(t,s,\xi),\qquad \mathcal Q(s,s,\xi)=I,
\end{equation}
which can be solved directly by the Peano-Baker series
\begin{equation}
   \mathcal Q(t,s,\xi) = I + \sum_{k=1}^\infty \i^k \int_s^t \mathcal R_2(t_1,s,\xi)\int_s^{t_1}\mathcal R_2(t_2,s,\xi)\cdots\int_s^{t_{k-1}}\mathcal R_2(t_k,s,\xi)\d t_k\cdots \d t_2\d t_1,
\end{equation}
where we used the notation $\mathcal R_2(t,s,\xi)=\tilde{\mathcal E}_0(s,t,\xi)R_2(t,\xi)\tilde {\mathcal E}_0(t,s,\xi)$. Note that $\|\mathcal R_2(t,s,\xi)\|=\|R_2(t,\xi)\|$ and therefore
\begin{equation}
\|\mathcal Q(t,s,\xi)\|\le \exp\big(\int_s^t\|R_2(\tau,\xi)\|\d\tau\big).
\end{equation}
Now the representation
\begin{equation}
   \mathcal M(t,\xi) = \frac{\lambda(t)}{\lambda(t+T)} M N_1(t+T,\xi) \tilde{\mathcal E}_0(t+T,t,\xi) \mathcal Q(t+T,t,\xi) N_1^{-1}(t,\xi)M^{-1}
\end{equation}
in combination with $\lambda(t)/\lambda(t+T)=\exp\beta T$ implies the desired result $\|\mathcal M(t,\xi)\|<1$ uniformly in $t\in[0,T]$ and $|\xi|\ge N$.

\subsection{Treatment of bounded frequencies}\label{sec1.2}
Our next aim is to show that for any $c>0$ there exists a natural number $k$ such that 
\begin{equation}
  \sup_{t\in[0,T]} \|\mathcal M^k(t,\xi)\| < 1
\end{equation}
for all $c\le|\xi|\le N$. Note, that we are only interested in a compact set in $t$ and $\xi$ here. We will 
combine spectral theory with a compactness argument. First step is to show that the spectrum $\spec\mathcal M(t,\xi)$ is contained inside the open unit ball $\{\zeta\in\C\,|\,|\zeta|<1\}$.

Note first, that $\mathcal M(t,\xi)\mathcal E(t,0,\xi)=\mathcal E(t+T,0,\xi) =\mathcal E(t,0,\xi) \mathcal M(0,\xi)$ and therefore $\mathcal M(t,\xi)$ is similar to $\mathcal M(0,\xi)$. Especially the spectrum
$\spec \mathcal M(t,\xi)$ is independent of $t$. Furthermore, we assumed $b(t)$ to be real. 
Therefore  \eqref{eq:ODE} has real-valued solutions and it follows that  $\mathcal M(t,\xi)$ is similar to a real matrix. Furthermore, according to Liouville theorem we see that 
\begin{equation}
    \det\mathcal M(t,\xi) = \exp\big(\i\int_t^{t+T} \trace A(\tau,\xi)\d\tau \big) = \exp\big(-2\int_t^{t+T} b(\tau)\d\tau\big) = \exp(-2\beta T),
\end{equation}
such that the eigenvalues $\varkappa_1(\xi)$ and $\varkappa_2(\xi)$ of 
$\mathcal M(t,\xi)$ are either real and of the form $\varkappa_2(\xi)=\varkappa_1(\xi)^{-1}\exp(-2\beta T)$ or complex-conjugate $\varkappa_2(\xi)=\overline{\varkappa_1(\xi)}$ and therefore $|\varkappa_1(\xi)|=|\varkappa_2(\xi)|=\exp(-\beta T)$. 
For the complex case we are done and $\spec\mathcal M(t,\xi) \subseteq \{ |\zeta|=\exp(-\beta T)\}$; for the real case we have to look more carefully. Note, that the eigenvalues  are continuous in $\xi$.

Assume that for a certain frequency $\bar\xi\ne0$ the monodromy matrix $\mathcal M(0,\bar\xi)$ has an eigenvalue with modulus $1$. Since it must be real, it is either $1$ or $-1$. Let $\vec c=(c_1,c_2)$ be a corresponding eigenvector. Then we can find a domain $\Omega_R=\{ x\in\R^n\,|\,|x|\le R\}$ such that 
$-|\bar\xi|^2$ is an eigenvalue of the Dirichlet Laplacian on $\Omega_R$ with normalised eigenfunction 
$\phi(x)$. If we consider the initial boundary value problem $\square u+2b(t)u_t = 0$ with Dirichlet boundary condition and $u(0,\cdot)=c_1\phi(x)$ and $u_t(0,\cdot)=\i c_2\phi(x)$, the corresponding solution 
satisfies 
\begin{equation}
  \begin{pmatrix} |\bar\xi| u(t,x) \\ \D_t u(t,x) \end{pmatrix}\bigg|_{t=T} = \mathcal M(T,\bar\xi) \vec c \phi(x) = \pm \vec c \phi(x)
\end{equation}
(sign according to the eigenvalue) and thus can be written as
\begin{equation}
 u(t,x) = f(t) \phi(x)
\end{equation}
with a $T$-periodic (or $2T$-periodic) function $f(t)$ and $f(0)=c_1$. However, this is not possible. If we denote the energy of this solution as $E(u;t) = \|\nabla u(t,\cdot)\|_2^2 + \|\D_t u(t,\cdot)\|_2^2$
the standard integration by parts argument gives
\begin{equation}
   \frac{\d}{\d t} E(u;t) = -2b(t) \|u_t\|_2^2 = -2 b(t) |f'(t)|^2,
\end{equation} 
such that after integration $0 = -2 \int_0^{(2)T} b(t) |f'(t)|^2\d t$. The positivity of $b(t)$ implies that $f$ is constant and this contradicts $\bar\xi\ne0$. Thus, $\pm 1\not\in\spec\mathcal M(t,\xi)$ for $\xi\ne0$
and therefore the spectral radius satisfies $\rho(\mathcal M(t,\xi))<1$. Thus the spectral radius formula
$\|\mathcal M^k(t,\xi)\|^{1/k} \to \rho(\mathcal M(t,\xi))<1$ implies that for any $t$ and $\xi$ we find a number $k$ such that $\|\mathcal M^k(t,\xi)\|<1$. 

Next, we want to show that we can find such a number $k$ uniform on any compact frequency interval
$|\xi|\in[c,N]$. Set for this $\mathcal U_k = \{ (t,\xi)\,|\, \|\mathcal M^{(2^k)}(t,\xi)\|<1\}$. The sets $\mathcal U_k$ are clearly open (by the continuity of the monodromy matrix) and satisfy $\mathcal U_k\subseteq\mathcal U_\ell$ for $k\le\ell$. The above reasoning shows that the compact set
$\mathcal C = \{ (t,\xi) \,|\, 0\le t\le T,\; c\le |\xi|\le N\}$ is contained in $\bigcup_k\mathcal U_k$ and by compactness we find one $k$ such that $\mathcal C\subset \mathcal U_k$. Hence the assertion of this section is proved.

{\bf Remark.} We know even a little bit more about the structure of the eigenvalues $\varkappa_1(\xi)$ and $\varkappa_2(\xi)$. We can apply a Liouville type transform to equation \eqref{eq:ODE} to deduce
Hill's equation
\begin{equation}\label{eq:Hill}
  v_{tt} + (|\xi|^2-b^2(t)-b'(t)) v =0,\qquad  v = \lambda(t) \hat u,
\end{equation}
such that Floquet theory, see e.g. \cite{Eastham2}, may be applied. This implies that if 
$b^2(t)+b'(t)$ is not constant (which is equivalent to $b(t)$ not constant) 
then there exist infinitely many intervals $I_0=(-\infty,\tau_0]$ and $I_k=[\tau_k^-,\tau_k^+]$, $k=1,2,\dots$, such that for $|\xi|\in I_j$, $j=0,1,\dots$, the spectrum $\spec\mathcal M(t,\xi)$ is real ({\em intervals of
instability} for \eqref{eq:Hill}), while for all other $\xi$ the eigenvalues are complex and conjugate to each other ({\em intervals of stability} for \eqref{eq:Hill}). The numbers $\tau_k^\pm$ are the eigenvalues of the corresponding periodic eigenvalue problem
$-v'' +(b^2(t)+b'(t))v=\lambda^2 v$ with periodic boundary conditions $v(0)=v(2T)$, $v'(0)=v'(2T)$.

\subsection{The neighbourhood of $\xi=0$}\label{sec13}
 The frequency $\xi=0$ is the only exceptional point of our 
reasoning, since $\spec\mathcal M(t,0) = \{1,\exp(-2\beta T)\}$ contains the eigenvalue $1$. This follows
directly by solving \eqref{eq:ODE}; a fundamental system of solutions is given by $1$ and $\int_0^t  \exp(-2\int_0^s b(\tau)\d\tau)\d s$. We will use ideas from the theory of Hill's equation to understand the structure of $\mathcal E(t,s,\xi)$ near $\xi=0$. From the remarks of Section~\ref{sec1.2} it is clear that $0$ is interior point of $I_0$
and, chosing $c$ small enough allows to write down a fundamental system of solutions to \eqref{eq:ODE} as
\begin{equation}
   \e^{-\nu_\pm(\xi) t} f_\pm(t,\xi) 
\end{equation}
with $T$-periodic functions $f_\pm(t,\xi)$ and exponents $\nu_\pm(\xi)$ such that $\exp(-\nu_\pm(\xi)T)\in\spec\mathcal M(t,\xi)$. It is clear that $\nu_\pm(\xi)>0$ for $\xi\ne0$ and we denote them in such a way that $\nu_+(\xi)\to0$ and $\nu_-(\xi)\to2\beta$  as $\xi\to0$. Any solution to \eqref{eq:ODE} is a combination of these two solutions. The part corresponding to $\nu_-(\xi)$ is not of interest for us
(because it leads to an exponential decay as $t\to\infty$) and we can concentrate on the $\nu_+(\xi)$ part.

We know that $\nu_+(\xi)$ is an analytic function of $|\xi|$ (as long as $\mathcal M(t,\xi)$ has no multiple eigenvalues) and can therefore be expanded into a MacLaurin series
\begin{equation}
   \nu_+(\xi) = \sum_{k=1}^\infty \alpha_k |\xi|^{k}.
\end{equation}
We want to show that $\alpha_1=0$ and $\alpha_2>0$. For this we use $\mathcal M(0,\xi) = \mathcal E(T,0,\xi)$ and calculate the derivatives of $\trace \mathcal M(0,\xi)$ with respect to $|\xi|$ at $\xi=0$.
Note, that $\partial_{|\xi|} A(t,\xi) = \big( \begin{smallmatrix}0&1\\1&0\end{smallmatrix}\big)=\mathcal J$ and $\partial_|\xi|^2A(t,\xi)=0$, such that  $\D_t\partial_{|\xi|} \mathcal E(t,s,\xi)=
\mathcal J\mathcal E(t,s,\xi) + A(t,\xi)\partial_{|\xi|}\mathcal E(t,s,\xi)$ and $\partial_{|\xi|}\mathcal E(s,s,\xi)=0$. Therefore we obtain the representation 
\begin{equation}
   \partial_{|\xi|} \mathcal E(t,s,\xi) = \int_s^t \mathcal E(t,\tau,\xi) \mathcal J \mathcal E(\tau,s,\xi)\d\tau
\end{equation}
Using that
\begin{equation} 
\mathcal E(t,s,0) = \exp\left(\i \int_s^t A(\tau,0)\d\tau\right) = 
\diag\left(1,\exp(-2\int_s^t b(\tau)\d\tau\right) = \diag\left(1,\frac{\lambda^2(s)}{\lambda^2(t)}\right)
\end{equation}
is diagonal, we immediately see that the above integrand has zeros as diagonal entries. This implies $\partial_{|\xi|}\trace \mathcal M(t,0)=0$. For the second derivative we use in analogy that
$\D_t\partial_{|\xi|}^2 \mathcal E(t,s,\xi) = A(t,\xi)\partial_{|\xi|}^2\mathcal E(t,s,\xi) + 2\mathcal J\partial_{|\xi|}\mathcal E(t,s,\xi)$ and $\partial_{|\xi|}^2\mathcal E(s,s,\xi)=0$, such that after integration
\begin{align}
  \partial_{|\xi|}^2 \mathcal E(t,s,\xi) &=2 \int_s^t \mathcal E(t,\tau,\xi) \mathcal J \partial_{|\xi|} \mathcal E(\tau,s,\xi)\d\tau \\
  &= 2\int_s^t \mathcal E(t,\tau,\xi) \mathcal J  \int_s^\tau \mathcal E(\tau,\theta,\xi) \mathcal J \mathcal E(\theta,s,\xi)\d\theta\d\tau\notag.
\end{align}
For $\xi=0$ we can evaluate these integrals and obtain for the trace 
\begin{equation}
  \partial_{|\xi|}^2\trace\mathcal M(t,0) = 2 \int_0^T\int_0^\tau \left(\frac{\lambda^2(\theta)}{\lambda^2(\tau)}+\frac{\lambda^2(\tau)}{\lambda^2(T)\lambda^2(\theta)}\right) \d\theta\d\tau>0.
\end{equation}
On the other hand, $\trace \mathcal M(t,\xi) = \exp(-\nu_+(\xi)T)+\exp(-\nu_-(\xi)T)$ with 
$\mathcal \nu_+(\xi)+\mathcal\nu_-(\xi) = 2\beta$, such that
$\partial_{|\xi|}\trace\mathcal M(t,0) = \alpha_1 T (1-\exp(-2\beta T)) = 0$ implies $\alpha_1=0$
and $\partial_{|\xi|}^2\trace \mathcal M(t,0) = 2\alpha_2 T (1-\exp(-2\beta T))>0$ implies $\alpha_2>0$.

Hence, we have shown that as $\xi\to0$ the exponent behaves like $\nu_+(\xi) = \alpha_2|\xi|^2+\mathcal O(|\xi|^3)$ (and, if we look carefully at the representations, we see that all odd coefficients vanish and thus the remainder term is $\mathcal O(|\xi|^4)$).
This will be enough to obtain energy and dispersive estimates for solutions to our Cauchy problem
in Section~\ref{sec2}. 

\subsection{Collection of results} What have we obtained so far? The main results are concerned with
the monodromy operator $\mathcal M(t,\xi) = \mathcal E(t+T,t,\xi)$ and its spectral properties.

\begin{thm}\label{thm1}
\begin{enumerate}
\item There exists a (large) number $N>0$ such that for all $|\xi|\ge N$ the monodromy matrix
$\mathcal M(t,\xi)$ is a contraction (uniform in $t$), i.e. 
$$\sup_t \|\mathcal M(t,\xi)\|<1.$$
\item For any  (small) number $c>0$ there exists
an exponent $k\in\mathbb N$, such that for $c\le |\xi|\le N$ the matrix $\mathcal M^k(t,\xi)$ is a contraction (uniform in $t$), i.e.
$\sup_t \|\mathcal M^k(t,\xi)\|<1$.
\item As $\xi\to0$ the eigenvalues of $\mathcal M(t,\xi)$ satisfy 
$$\log \varkappa_1(\xi) = -\alpha_2T|\xi|^2 + \mathcal O(|\xi|^4), \qquad
\log \varkappa_2(\xi) = -2\beta T+\alpha_2T|\xi|^2 + \mathcal O(|\xi|^4) $$
with a positive coefficient $\alpha_2>0$.
\end{enumerate}
\end{thm}

\section{Estimates for solutions}\label{sec2}
The repesentations from Section~\ref{sec1} allow us to estimate the Fourier transform of solutions,
in combination with Plancherel's theorem this gives estimates in $L^2$-spaces, combined with H\"older
inequality and mapping properties of the Fourier transform dispersive estimates follow.

\subsection{Energy estimates}  
We distinguish between small and large frequencies. If $|\xi|\ge N$ the monodromy matrix $\mathcal M(0,\xi)$ is a contraction and therefore $\|\mathcal E(t,0,\xi)\|=\|\mathcal M^\ell(s,\xi)\mathcal E(s,0,\xi)\|\le c^\ell \|\mathcal E(s,0,\xi)\|$ for $t=\ell T+s$, $s\in[0,T]$ and $c=\sup_t\|\mathcal M(t,\xi)\|<1$. Furthermore, since
$b(t)\ge 0$ we know that $\|\mathcal E(s,0,\xi)\|\le 1$ and therefore
\begin{equation}
   \|\mathcal E(t,0,\xi)\| \le \e^{-\delta (t-T)}
\end{equation}
with $\delta = T^{-1} \log c^{-1}>0$. Thus, high frequencies lead to an exponential decay. For the intermediate frequencies we obtain similarly  $\|\mathcal E(t,0,\xi)\|=\|\mathcal M^\ell(s,\xi)\mathcal E(s,0,\xi)\|\le c^\ell \|\mathcal E(s,0,\xi)\|$ for $t=\ell k T+s$, $s\in[0,kT]$ and $c=\sup_t\|\mathcal M(t,\xi)\|<1$.
Again this yields exponential decay, but now of the form
\begin{equation}
   \|\mathcal E(t,0,\xi)\| \le \e^{-\delta (t-kT)}
\end{equation}
with $\delta = (kT)^{-1} \log c^{-1}>0$. Hence, the only non-exponential contribution may come from
the neighbourhood of $\xi=0$. For the treatment of small frequencies we have to specify the structure of the estimate we have in mind. While estimating the energy of the solution at time $t$ in terms of the initial energy brings (due to $1\in\spec\mathcal M(t,0)$) only the trivial uniform bound and no decay, an estimate in terms of $\|u_1\|_{H^1}$ and $\|u_2\|_{L^2}$ brings decay. Reason for this is that we can use an additional factor $|\xi|$ for small frequencies.

If $|\xi|\le c$ is sufficiently small we know that a fundamental system of solutions to \eqref{eq:ODE} is given by $\exp(-\nu_\pm(\xi)t)f_\pm(t,\xi)$ with $T$-periodic functions $f_\pm(t,\xi)$ and exponents $\nu_+(\xi)=\alpha_2|\xi|^2+\mathcal O(|\xi|^4)$, $\nu_-(\xi)=2\beta-\alpha_2|\xi|^2+\mathcal O(|\xi|^4)$. Furthermore,
$f_\pm(t,\xi)$ are non-zero for all $t$ and $\xi$. This follows from the fact that they are periodic and non-zero for $\xi=0$. Thus, if they would have a zero for some $t$ and $\xi$ we could find a smallest value of $\xi$ where the zero occurs. By differentiability of $f_\pm$ it follows that for this fixed $\xi$ we would obtain a zero of order at least 2, which contradicts the fact that $\exp(-\nu_\pm(\xi)t)f_\pm(t,\xi)$ is a not identically vanishing solution of the second order equation \eqref{eq:ODE}. Hence, we may assume that $f_\pm(0,\xi)=1$ for all $|\xi|\le c$. This allows to express the special fundamental system of solutions
$\Phi_1(t,\xi)$ and $\Phi_2(t,\xi)$ with $\Phi_1(0,\xi)=1$, $\partial_t\Phi_1(0,\xi)=0$,
$\Phi_2(0,\xi)=0$ and $\partial_t\Phi_2(t,\xi)=1$, i.e. the fundamental system representing solutions
to \eqref{eq:ODE} as
\begin{equation}
  \hat u(t,\xi) = \sum_{j=1,2} \Phi_j(t,\xi)\hat u_j(\xi),\qquad |\xi|\le c,
\end{equation} 
in terms of $f_\pm(t,\xi)$ and the exponents $\nu_\pm(\xi)$. A simple calculation shows
\begin{align}
\Phi_2(t,\xi)&= \frac{\e^{-\nu_+(\xi)t}f_+(t,\xi)-\e^{-\nu_-(\xi)t}f_-(t,\xi)}{\nu_-(\xi)-\nu_+(\xi) + \partial_tf_+(0,\xi)-\partial_t f_-(0,\xi)},\\
\Phi_1(t,\xi) &=\frac{\e^{-\nu_+(\xi)t}f_+(t,\xi)+\e^{-\nu_-(\xi)t}f_-(t,\xi)}2-
\bigg(\frac{\partial_tf_+(0,\xi)+\partial_tf_-(0,\xi)}2-\beta\bigg)\Phi_2(t,\xi).
\end{align}
Both functions are smooth in $\xi$ and differentiable in $t$, especially it follows that the denominator in the first expression is non-zero. Since we are interested in polynomial decay rates, we can forget about all the $f_-$-terms, which immediately lead to exponential decay. Thus, to estimate $|\xi|\hat u(t,\xi)$
in terms of $\hat u_1$ and $\hat u_2$, the typical term to estimate is $|\xi|\e^{-\nu_+(\xi)t}f_+(t,\xi)$
(multiplied by  a $\xi$-dependent function). Now $\nu_+(\xi)\sim\alpha_2|\xi|^2$ implies the uniform
decay rate $t^{-1/2}$ for this term. To estimate $\partial_t \hat u(t,\xi)$ in terms of $\hat u_1$ and $\hat u_2$ we have to consider the typical term 
$\e^{-\nu_+(\xi)t} \partial_t f_+(t,\xi) - \nu_+(\xi) \e^{-\nu_+(\xi)t}f_+(t,\xi)$. The second addend gives $t^{-1}$, while the first one has to be considered in detail.  Note, that $\Phi_1(t,0)=1$, such that
comparing representations implies $f_+(t,0)=1$.
Since the equation was parametrised by $|\xi|^2$ smoothness in $|\xi|$ and periodicity in $t$ imply $\partial_t f_+(t,\xi)=|\xi|^2\tilde h_+(t,\xi)$ with a 
bounded $T$-periodic function $\tilde h_+(t,\xi)$. Therefore we see, that the first addend gives the same decay rate $t^{-1}$. 

We collect our results in the following theorem.
\begin{thm}\label{thm2}
The solution $u(t,x)$ of the Cauchy problem \eqref{eq:CP} satisfies the a-priori estimates
\begin{align*}
 \|u(t,\cdot)\|_{L^2}&\lesssim \|u_1\|_{L^2} + \|u_2\|_{H^{-1}}, \\
 \|\nabla u(t,\cdot)\|_{L^2} &\lesssim (1+t)^{-1/2} \big(\|u_1\|_{H^1} + \|u_2\|_{L^2} \big),\\
  \|\partial_t u(t,\cdot)\|_{L^2} &\lesssim (1+t)^{-1} \big(\|u_1\|_{H^1} + \|u_2\|_{L^2} \big).
\end{align*}
Furthermore, for any cut-off function $\chi\in C^\infty(\R)$ with $\chi(s)=0$ near $s=0$ and
$\chi(s)=1$ for large $s$ exists a constant $\delta>0$ such that the exponential estimate
$$ \|\chi(\D) u(t,\cdot)\|_{L^2} + \|\chi(\D) \nabla u(t,\cdot)\|_{L^2} + \|\chi(\D)\partial_t u(t,\cdot)\|_{L^2}
\lesssim \e^{-\delta t}  \big(\|u_1\|_{H^1} + \|u_2\|_{L^2} \big) $$
holds true.
\end{thm}

\subsection{Dispersive estimates} We will continue this short note with some remarks on dispersive
and more generally $L^p$--$L^q$ decay estimates. Again only the small frequencies are of interest, since by Sobolev embedding the previous theorem implies
\begin{equation}
 \|\chi(\D) u(t,\cdot)\|_{L^q} + \|\chi(\D) \nabla u(t,\cdot)\|_{L^q} + \|\chi(\D)\partial_t u(t,\cdot)\|_{L^q}
\lesssim \e^{-\delta t}  \big(\|u_1\|_{H^{p,r_p+1}} + \|u_2\|_{H^{p,r_p}} \big) 
\end{equation}
for any choice of indices $1\le p\le 2\le q\le \infty$ and regularity $r_p>n(1/p-1/q)$. Thus it remains to
consider the typical terms from the previous section near $\xi=0$. Instead of Plancherel's theorem
we use H\"older inequality together with the $L^p$--$L^{p'}$ boundedness of the Fourier transform
for $pp'=p+p'$ to deduce
\begin{equation}
  \|\chi(\D) |\D| \e^{-\nu_+(\D)} f_+(t,\D) \|_{L^p\to L^q}
  \le \| \chi(\xi) |\xi| \e^{-\nu_+(\xi)}f_+(t,\xi) \|_{L^r}
\end{equation}
for any $1\le p\le 2\le q\le \infty$ and with $\frac1r = \frac1p-\frac1q$. The $L^r$-norm can be calculated
directly using $\nu_+(\xi)\sim\alpha_2|\xi|^2$, which implies the decay rate $t^{-1/2-n/2r}$. Similarly we obtain for the derivative terms the rate $t^{-1-n/2r}$ and for the solution itself $t^{-n/2r}$.

\begin{thm}\label{thm3}
The solution $u(t,x)$ of the Cauchy problem \eqref{eq:CP} satisfies the a-priori estimates
\begin{align*}
 \|u(t,\cdot)\|_{L^q}&\lesssim (1+t)^{-\frac n2(\frac1p-\frac1q)} \big(\|u_1\|_{H^{p,r_p}} + \|u_2\|_{H^{p,r_p-1}}\big), \\
 \|\nabla u(t,\cdot)\|_{L^q} &\lesssim (1+t)^{-\frac12-\frac n2(\frac1p-\frac1q)} \big(\|u_1\|_{H^{p,r_p+1}} + \|u_2\|_{H^{p,r_p}} \big),\\
  \|\partial_t u(t,\cdot)\|_{L^q} &\lesssim (1+t)^{-1-\frac n2(\frac1p-\frac1q)} \big(\|u_1\|_{H^{p,r_p+1}} + \|u_2\|_{H^{p,r_p}} \big).
\end{align*}
for all $1\le p\le 2\le q\le \infty$ and $r_p>n(1/p-1/q)$.
\end{thm}

\subsection{Diffusion phenomenon} For proving estimates for the solution $u$ we used that
the only bad term in the representation of solutions was $\e^{-\nu_+(\xi)t} f_+(t,\xi)
\sim \e^{-\alpha_2|\xi|^2t}$, which corresponds to the Fourier multiplier for a corresponding
heat equation 
\begin{equation}\label{eq:heat}
   w_t = \alpha_2 \Delta w, \qquad w(0,\cdot)=w_0.
\end{equation}
Choosing $w_0$ in dependence of $u_1$ and $u_2$ allows to cancel the corresponding terms
in the representation of solutions such that the norm of the difference $||u-w||_{L^2}$ decays. For
constant $b$ and with $\alpha_2=(2b)^{-1}$ this was observed in \cite{Nis97} and \cite{HM01}. For periodic
dissipation terms a similar statement is valid.
If we choose
\begin{equation}\label{eq:w0}
   w_0 = \left(\frac12+\frac{\beta-\gamma}{2\beta-\gamma}\right)u_1 
   + \frac1{2\beta -\gamma} u_2,\qquad \gamma=\partial_tf_-(0,0)=2\beta-(1-\e^{-2\beta T})(\int_0^T \frac{\d\tau}{\lambda^2(\tau)})^{-1}
\end{equation}
and use the $\alpha_2$ from Section~\ref{sec13},
\begin{equation}
\alpha_2 = \frac1{T(1-\e^{-2\beta T})} \int_0^T \int_0^\tau\left( \frac{\lambda^2(\theta)}{\lambda^2(\tau)}+\frac{\lambda^2(\tau)}{\lambda^2(T)\lambda^2(\theta)}\right)\d\theta\d\tau,
\qquad \lambda(t)=\exp\left(\int_0^t b(s)\d s\right),
\end{equation}
then the following result holds true:
\begin{thm}\label{thm4}
   The solutions $u(t,x)$ of \eqref{eq:CP} and $w(t,x)$ of \eqref{eq:heat} satisfy under 
   the relation \eqref{eq:w0} the a-priori estimate
   \begin{equation}
      ||u(t,\cdot)-w(t,\cdot)||_{L^2} \lesssim (1+t)^{-1} \big(|| u_1||_{H^1}+||u_2||_{L^2}\big).
   \end{equation}
 \end{thm}

To prove this result we first note that we can forget about all terms in the representation which
give a faster decay. The choice of the initial datum \eqref{eq:w0} implies that the only term of interest
is $(\e^{-\nu_+(\xi)t}f_+(t,\xi) - \e^{-\alpha_2|\xi|^2t})\hat w_0$ and $f_+(t,\xi)=1+\mathcal O(|\xi|^2)$
together with $\nu_+(\xi)=\alpha_2|\xi|^2+\mathcal O(|\xi|^4)$ localised near $|\xi|=0$. But this multiplier can be estimated by
a combination of $\e^{-\nu_+(\xi)t} |\xi|^2$ and $\e^{-\min(\nu_+(\xi),\alpha_2|\xi|^2)t}|\xi|^4 t$. Both terms 
decay uniformly like $(1+t)^{-1}$ and the assertion follows.

A similar statement with improvement of one decay order holds for dispersive and $L^p$--$L^q$ estimates as well as for estimates of higher order spatial derivatives. The reasoning is analogous.

\bigskip
{\bf Acknowledgements.} Ideas to this note and several arguments have been achieved while the author visited Prof. Ryo Ikehata at Hiroshima Unicersity financed by a joint DFG / DMV grant 2005. The author is also grateful to Prof. Kenji Nishihara for many discussions and valuable hints on the diffusion phenomenon for damped waves.


\begin{thebibliography}{99}
\bibitem[Eas73]{Eastham2}
M.S.P. Eastham, {\it The spectral theory of periodic differential equations},
Scottish Academy Press, Edinburgh and London, 1973.
\bibitem[Eas89]{Eastham}
\bysame, {\it The asymptotic solution of linear differential equations},
Oxford Science Publications, 1989.
\bibitem[HiRe03]{ReissigHirosawa}
F. Hirosawa, M. Reissig, {\em From wave to Klein-Gordon type decay rates}
in: {\em Nonlinear hyperbolic equations, spectral theory and wavelet transformations} (S. Albeverio, M. Demuth, E. Schrohe, B.-W. Schulze ed.), Operator Theory, Advances and Applications, Vol. 145, p. 95--155, Birkh\"auser Verlag, Basel 2003
\bibitem[Mat76]{Matsumura76}
A. Matsumura, {\em On the asymptotic behavior of solutions of
dissipative wave equations}, Publ. Res. Inst. Math. Sci. 12/1 (1976) 169--189.
\bibitem[Nis97]{Nis97} 
K. Nishihara,  {\em Asymptotic behaviour of solutions of quasilinear hyperbolic equations  with linear
damping}, J. Differential Equations 137/2 (1997) 384--395.
\bibitem[Rei04]{Reissig}
M. Reissig, {\em $L_p$--$L_q$ decay estimates for wave equations with time-dependent coefficients}, J. Nonlin. Math. Phys. 11/4 (2004) 534--548.
\bibitem[ReSm05]{ReissigSmith}
M. Reissig, J. Smith, {\em $L^p$--$L^q$ estimate for wave equation
with bounded time-dependent coefficient}, Hokkaido Math. J. 34/3 (2005) 541--586.
\bibitem[ReYa00]{ReissigYagdjian}
M. Reissig, K. Yagdjian, {\em About the influence of oscillations on Strichartz type decay estimates}, Rend. Sem. Mat. Univ. Pol. Torino, 58/3 (2000) 375--388.
\bibitem[Wir06]{Wirth06}
J. Wirth, {\it Wave equations with time-dependent dissipation I. Non-effective dissipation}, J. Differential Equations, 222/2 (2006) 487--514.
\bibitem[Wir07]{Wirth07}
\bysame, {\it Wave equations with time-dependent dissipation II. Effective dissipation}, J. Differential Equations, 232/1 (2007) 74--103.
\bibitem[Yag01] {Yagdjian01}
K. Yagdjian, {\em Parametric resonance and nonexistence of the global solution to nonlinear wave equations}, J. Math. Anal. Appl. 260/1 (2001) 251--268.
\bibitem[YaMi01]{HM01} 
Han Yang, A. Milani, {\em On the diffusion phenomenon of quasilinear hyperbolic waves}, Bull. Sci. Math. 124/5 (2000) 415--433.
\end{thebibliography}
\end{document}